\documentclass[12pt]{amsart}

\usepackage{amsmath,amsthm,amssymb}
\usepackage{tikz}
\usetikzlibrary{arrows.meta}

\usepackage{a4wide}
\usepackage{enumerate}
\usepackage{multicol, multirow}
\newtheorem{theorem}{Theorem}[section]
\newtheorem{lemma}[theorem]{Lemma}

\newtheorem{theo}{Theorem}

\theoremstyle{remark}
\newtheorem{remark}[theorem]{Remark}

\newtheorem*{claim*}{Claim}

\newcommand{\C}{\ensuremath{\mathbb{C}}}
\newcommand{\R}{\ensuremath{\mathbb{R}}}

\newcommand{\g}[1]{\ensuremath{\mathfrak{#1}}}

\newcommand{\cal}[1]{\ensuremath{\mathcal{#1}}}

\DeclareMathOperator{\tr}{tr}

\DeclareMathOperator{\Ad}{Ad}
\DeclareMathOperator{\Aut}{Aut}

\DeclareMathOperator{\Exp}{Exp}

\DeclareMathOperator{\spann}{span}

\DeclareMathOperator{\rank}{rank}

\begin{document}
\title[Isoparametric hypersurfaces in symmetric spaces of higher rank]{Isoparametric hypersurfaces in symmetric spaces\\ of non-compact type and higher rank}
\author[M.~Dom\'{\i}nguez-V\'{a}zquez]{Miguel Dom\'{\i}nguez-V\'{a}zquez}
\address{Department of Mathematics, Universidade de Santiago de Compostela, R\'ua Lope G\'omez de Marzoa s/n, 15782 Santiago de Compostela, Spain.
\newline\indent 
CITMAga, 15782 Santiago de Compostela, Spain.}
\email{miguel.dominguez@usc.es}
\author[V.~Sanmart\'in-L\'opez]{V\'ictor Sanmart\'in-L\'opez}
\address{Department of Mathematics, Universidade de Santiago de Compostela, R\'ua Lope G\'omez de Marzoa s/n, 15782 Santiago de Compostela, Spain.}
\email{victor.sanmartin@usc.es}

\begin{abstract}
We construct inhomogeneous isoparametric families of hypersurfaces with non-austere focal set on each symmetric space of non-compact type and rank $\geq 3$. If the rank is $\geq 4$, there are infinitely many such examples. Our construction yields the first examples of isoparametric families on any Riemannian manifold known to have a non-austere focal set. They can be obtained from a new general extension method of submanifolds from Euclidean spaces to symmetric spaces of non-compact type. This method preserves the mean curvature and isoparametricity, among other geometric properties.
\end{abstract}

\thanks{The authors have been supported by projects PID2019-105138GB-C21, PID2022-138988NB-I00 (MCIN/AEI/10.13039/501100011033/FEDER, EU, Spain) and ED431F 2020/04 (Xunta de Galicia, Spain). The first author acknowledges support of the Ram\'{o}n y Cajal program of the Spanish State Research Agency. 
} 

\subjclass[2020]{53C40, 53C35, 53C42}
\keywords{Isoparametric hypersurface, inhomogeneous, austere submanifold, symmetric space, non-compact type, hyperpolar, extension}
\maketitle
\vspace{-0.5ex}
\section{Introduction}

An isoparametric family of hypersurfaces is a decomposition of a Riemannian manifold into equidistant hypersurfaces of constant mean curvature and possibly one or two focal submanifolds with codimension greater than one. An isoparametric family is called homogeneous if it is the orbit foliation of an isometric action.
In this article we provide a construction method of isoparametric hypersurfaces from where the next result follows.

\begin{theo}\label{th:first}
Each symmetric space $M$ of non-compact type and rank at least three admits inhomogeneous isoparametric families of hypersurfaces with non-austere focal submanifolds. If the rank is greater than or equal to four, there exist uncountably many such examples, up to congruence.
\end{theo}

We emphasize three novel and surprising features of this result. 

Firstly, among the symmetric spaces of non-compact type, only those in the very special class of  $(BC_r)$ type (namely, those with non-reduced root system) were known to admit inhomogeneous isoparametric hypersurfaces~\cite{DRDV:adv,DV:imrn,DDR:hhn}. In this paper we construct examples in symmetric spaces of any type and rank at least three. 

Secondly, Theorem~\ref{th:first} provides the first isoparametric hypersurfaces which are known to have a non-austere focal submanifold in any Riemannian manifold. We recall that a submanifold is said to be austere if its multiset of principal curvatures is invariant under change of sign~\cite{HL}. In particular, austere submanifolds are minimal. The focal submanifolds of an isoparametric family of hypersurfaces are always minimal~\cite{GT:asian}, and Theorem~\ref{th:first} shows that they do not need to be austere. 
We highlight the contrast with homogeneous hypersurfaces, as singular orbits of cohomogeneity one actions are always austere~\cite{Podesta}. More generally, the focal submanifolds of isoparametric families of hypersurfaces with constant principal curvatures~\cite{GT:asian} (such as the famous inhomogeneous examples in spheres~\cite{FKM:mathz}), or of previous constructions in symmetric spaces of non-compact type~\cite{DRDV:adv,DV:imrn,DDR:hhn}, are known to be austere. Indeed, the austerity of the focal submanifolds has played a crucial role in the classification of isoparametric hypersurfaces in spheres~\cite{Mi:ann} and complex hyperbolic spaces~\cite{DDS:adv}. Theorem~\ref{th:first} reveals that this property does not hold in general, even in spaces with a large isometry group. For more information on isoparametric hypersurfaces, see~\cite{BCO,CR:book,Chi,Chi:survey}.

Thirdly, Theorem~\ref{th:first} also holds for reducible spaces and yields the  first inhomogeneous examples that are not extrinsic products in any reducible symmetric space. For instance, any product of at least three rank one symmetric spaces of non-compact type turns out to admit inhomogeneous isoparametric hypersurfaces. Particularly remarkable is the case of the product of three real hyperbolic planes, $\R H^2\times \R H^2\times \R H^2$, which admits such an inhomogeneous example, despite the fact that the isoparametric hypersurfaces of the irreducible factors are all homogeneous. Together with the complex hyperbolic space $\C H^3$~\cite{DRDV:adv}, this reducible space seems to be the lowest dimensional symmetric space that is known to admit inhomogeneous isoparametric hypersurfaces.

\medskip
We now describe the construction of the new examples. Let $M\cong G/K$ be a symmetric space of non-compact type, where $G$ is the connected component of the identity of the isometry group of $M$, and $K$ is the isotropy group at some point $o\in M$. Let $\g{g}=\g{k}\oplus\g{p}$ be the corresponding Cartan decomposition, and $\g{a}$ a maximal abelian subspace of $\g{p}$. Consider the restricted root space decomposition $\g{g}=\g{g}_0\oplus(\bigoplus_{\lambda\in\Delta}\g{g}_\lambda)$ with respect to $\g{a}$, where $\Delta\subset\g{a}^*$ is the set of restricted roots.   Let $\Delta^+$ be a choice of positive roots in $\Delta$, and define $\g{n}=\bigoplus_{\lambda\in\Delta^+}\g{g}_\lambda$. Consider the connected Lie subgroups $A$, $N$ and $AN$ of $G$ with Lie algebras $\g{a}$, $\g{n}$ and $\g{a}\oplus\g{n}$, respectively. The Iwasawa decomposition theorem for $G$ implies that $M\cong G/K$ is isometric to the solvable Lie group $AN$ equipped with a left-invariant metric $\langle\cdot,\cdot \rangle$. For each $\lambda\in \Delta$, we define $H_\lambda\in\g{a}$ by $\langle H,H_\lambda\rangle =\lambda(H)$ for all $H\in\g{a}$, and put
\[
H_\delta=\frac{1}{2}\sum_{\lambda\in\Delta^+}\dim\g{g}_\lambda\, H_\lambda.
\]
 
\begin{theo}\label{th:explicit}
	Let $M$ be a symmetric space of non-compact type and rank $\geq 3$. Let $S$ be the connected Lie subgroup of $AN$ with Lie algebra $\g{s}=\g{b}\oplus\g{n}$, where $\g{b}$ is any subspace of codimension at least two of $\g{a}$ such that $H_\delta\in\g{b}$. Then:
	\begin{enumerate}[\rm(i)]
		\item The orbit $S\cdot o$ is a minimal submanifold. It is non-austere for a generic choice of $\g{b}$ as above, or if $\dim \g{b}=1$.\
		\item The distance tubes around $S\cdot o$ define an inhomogeneous isoparametric family of hypersurfaces with non-constant principal curvatures on $M$.
		\item If two choices of subspaces $\g{b}_1$, $\g{b}_2$ of $\g{a}$ as above produce congruent isoparametric families, then $\g{b}_2=\varphi(\g{b}_1)$ for some  $\varphi$ in the finite group $\Aut(\Delta)$ of linear isometries of $\g{a}$ that leave the subset $\{H_\lambda\in\g{a}: \lambda\in\Delta\}$ invariant.
	\end{enumerate}
\end{theo}
The condition $H_\delta\in\g{b}\subset\g{a}$ ensures that the orbit $S\cdot o$ is not only minimal, but an Einstein solvmanifold whenever $M$ is Einstein (for instance, if $M$ is irreducible); see Remark~\ref{rem:Einstein}. Theorem~\ref{th:explicit} guarantees that any choice of subspace $\g{b}$ of codimension at least two in $\g{a}$ containing $H_\delta$ gives rise to an inhomogeneous isoparametric family of $M$.
The genericity condition in item (i) is only needed for the non-austerity of the focal submanifold ${S\cdot o}$. Here, by a generic choice we mean that $\g{b}$ can be selected in an open dense subset of the Grassmannian of $k$-dimensional subspaces of $\g{a}$ containing $H_\delta$, for some $k\in\{1,\dots,\rank M-2\}$. However, genericity can be substituted by explicit computable conditions in terms of the roots. In any case, it is easy to produce concrete examples of non-austere focal sets, the simplest method being just taking $\g{b}=\R H_\delta$. This implies the first claim of Theorem~\ref{th:first}. However, Theorem~\ref{th:explicit} provides more examples of isoparametric families than Theorem~\ref{th:first}, as tubes around $S\cdot o$ are always isoparametric with non-constant principal curvatures, independently on whether $S \cdot o$ is austere or not, and $S\cdot o$ is actually austere for certain choices of $\g{b}$. 
See the discussion around~\eqref{eq:austere_characterization}-\eqref{eq:sufficient_non_austere}, as well as Remark~\ref{rem:austere}, for further details. 

If $\rank M=\dim\g{a}\geq 4$, then there are uncountably many subspaces $\g{b}$ of codimension at least two in $\g{a}$ with $H_\delta\in\g{b}$. By the generic property in Theorem~\ref{th:explicit}~(i), an uncountable number of such subspaces produce non-austere submanifolds of the form $S\cdot o$. Now, according to Theorem~\ref{th:explicit}~(iii), tubes around the latter give rise to uncountably many congruence classes of isoparametric families with non-austere focal set. This proves the infiniteness claim in Theorem~\ref{th:first}. In conclusion, Theorem~\ref{th:explicit} clearly implies Theorem~\ref{th:first}. %
\medskip

We would like to point out that Theorem~\ref{th:explicit}~(iii) and its proof may play an important role in the determination of the congruence classes of homogeneous hyperpolar foliations~\cite{BDT:jdg}, which is still an open problem. The relation to hyperpolar foliations stems from the fact that the group $S$ acts hyperpolarly on $M$ (see paragraph below and \S\ref{subsec:extension}). Our proof of Theorem~\ref{th:explicit}~(iii) is conceptual and case-free, and does not make use of the assumptions on the rank of $M$, the codimension of $\g{b}$, or $H_\delta\in\g{b}$. The particular case of $\g{b}$ of codimension one in $\g{a}$ has been completely characterized in~\cite[pp.~9--20]{BT:jdg} and~\cite{So:arxiv} in terms of the automorphisms of the Dynkin diagram of $M$.


We will provide two proofs of the claims of  isoparametricity and inhomogeneity in Theorem~\ref{th:explicit}~(ii). The first one is a direct approach based on the study of the extrinsic geometry of the examples via standard Jacobi field theory. The second one arises as a straightforward consequence of a new extension method of submanifolds and foliations in symmetric spaces from totally geodesic, flat submanifolds.  This method is of independent interest, and seems to have remained unnoticed, in spite of being an application of results from~\cite{BDT:jdg} and~\cite{DV:imrn}. Specifically, the action of a group $S$ as in Theorem~\ref{th:explicit} on a symmetric space $M$ happens to be free, hyperpolar, and with minimal orbits~\cite{BDT:jdg}. Then, a general extension result (for free, polar actions with minimal orbits on Riemannian manifolds) given in~\cite{DV:imrn} by the first author applies, thus allowing to extend submanifolds and singular Riemannian foliations from a section $\Sigma$ of the polar $S$-action to~$M$. Importantly, this extension procedure preserves various geometric properties such as minimality, polarity or isoparametricity. We refer to~\cite{DV:imrn} (see also~\cite{AB:book} and~\cite{HLO}) for the definitions of isoparametric submanifold of arbitrary codimension, singular Riemannian foliation, and  polar, hyperpolar or isoparametric foliation.

\begin{theo}\label{th:extension}
	Let $M$ be a symmetric space of non-compact type and rank $\geq 3$, and $S$ the connected Lie subgroup of $AN$ with Lie algebra $\g{s}=\g{b}\oplus\g{n}$, where $\g{b}$ is a subspace of codimension at least two of $\g{a}$ such that $H_\delta\in\g{b}$.
	Let $\Sigma\cong \R^{d}$ be a section of the hyperpolar $S$-action on $M$, where $d:=\rank M-\dim\g{b}$.

	Let $P$ and $\cal{F}$ be a connected submanifold and a singular Riemannian foliation of codimension $k$ in the Euclidean space $\Sigma\cong \R^{d}$, respectively.~Then
		\[
		S\cdot P:=\{s(p):s\in S,\, p\in P\} \quad \text{and} \quad S\cdot \cal{F}:=\{S\cdot L: L\in \cal{F}\}
		\]
		are a connected submanifold and a singular Riemannian foliation of codimension $k$ in the symmetric space $M$, respectively. Moreover:
		\begin{enumerate}[\rm (i)]
			\item If $P$ has parallel mean curvature, is minimal, has flat normal bundle, or is isoparametric as a submanifold of $\Sigma$, the same is true for $S\cdot P$ as a submanifold~of~$M$. However, $S\cdot P$ is never totally geodesic.
			\item If $\cal{F}$ is a polar (or equivalently, isoparametric) singular Riemannian foliation of $\Sigma$, then $S\cdot \cal{F}$ is an isoparametric and hyperpolar singular Riemannian foliation of $M$.
			\item $S\cdot P$ (resp.\ $S\cdot\cal{F}$) is inhomogeneous in $M$, except if $P$ is an affine subspace (resp.\ except if $\cal{F}$ is a foliation given by parallel affine subspaces)  of $\Sigma\cong\R^d$.
		\end{enumerate}
\end{theo}

Under the perspective of the extension procedure described in Theorem~\ref{th:extension}, the isoparametric families of hypersurfaces constructed in Theorem~\ref{th:explicit} are nothing but extensions $S\cdot \cal{F}$ of the most basic examples of isoparametric foliations $\cal{F}$ with a singular leaf on a Euclidean space~$\R^d$, namely those given by concentric spheres and their common center. In other words, a tube of radius $t>0$ around $S\cdot o$ (which is a leaf of the isoparametric foliation mentioned in Theorem~\ref{th:explicit}~(ii)) is the union of $S$-orbits passing through the points of a sphere of radius $t$ in $\Sigma=\Exp(\g{b}^\perp)\cdot o$ centered at~$o$, where $\g{b}^\perp$ is the orthogonal complement of $\g{b}$ in $\g{a}$.

The general extension result in~\cite{DV:imrn} was applied, in that paper, to a different class of free polar actions with minimal orbits on symmetric spaces of non-compact type. 
The resulting extension method for symmetric spaces allowed to enlarge submanifolds and foliations from semisimple symmetric spaces of lower rank embedded in a particular totally geodesic manner (namely, as so-called boundary components). This contrasts with the method provided in Theorem~\ref{th:extension}, where one extends submanifolds or foliations from Euclidean spaces. Even more striking is the fact that, whereas the method in~\cite{DV:imrn} preserved the homogeneity of the examples, Theorem~\ref{th:extension}~(iii) reveals that this is almost never the case for the extension procedure we present in this article. The ultimate reason for this difference is that the connected component of the identity of the isometry group of a Euclidean space, embedded as a maximal flat in $M$, does not embed into the isometry group of $M$ (only translations do), whereas there does exist such an embedding for boundary components~of~$M$. 

At the topological level, the extension method in Theorem~\ref{th:extension} works as extending submanifolds of a subspace of a Euclidean space by perpendicularly attaching parallel affine subspaces, thus giving rise to a cylinder over the original submanifold. This is so because the $S$-orbit foliation on $M$ is diffeomorphically equivalent to a foliation of $\R^{\dim M}$ by affine $d$-codimensional subspaces. However, geometrically speaking, the extension procedure in Theorem~\ref{th:extension} is not trivial at all, in that the attached fibers (i.e.\ the $S$-orbits) are minimal, but never totally geodesic.

Although it is not the purpose of this paper to exploit the whole potential of Theorem~C to produce interesting examples of submanifolds in symmetric spaces of non-compact type, we will briefly illustrate it with some simple applications. For instance, by extending a catenoid in a $3$-dimensional Euclidean space, we can obtain a complete minimal hypersurface diffeomorphic to $\mathbb{S}^1\times \R^{\dim M-2}$ on any symmetric space of non-compact type and rank at least $4$. This can of course be generalized in several ways, for example by extending (extrinsic products of) minimal surfaces with more complicated topologies. In a different direction, by extending polar foliations on Euclidean spaces, one can obtain (to our knowledge, the first known examples of) inhomogeneous, isoparametric, hyperpolar, singular Riemannian foliations of codimension higher than one on a non-compact symmetric space $M$. Such foliations have no totally geodesic leaves, and can take any codimension $k\in\{1,\dots, \rank M-2\}$. Moreover, by considering a symmetric space $M$ of sufficiently large rank, one can construct examples with any prescribed polar infinitesimal foliation.

\medskip
\textbf{Acknowledgements.} We would like to thank J.~Berndt, E.~Garc\'ia-R\'io, J.~M.~Manzano, H.~Tamaru, and the anonymous referees for helpful comments and suggestions. We are grateful to I.~Solonenko for pointing out an issue in the study of the congruence in a previous version of this article.

\section{Proof of the main theorems}

We start by introducing some notation and known facts on symmetric spaces of non-compact type. For more information, we refer the reader to~\cite{DDS:sp} and \cite[Chapter~2]{eberlein}, and to~\cite[Sections~II.5, VI.4-5]{Knapp} for the theory of root systems.

As in the introduction, let $M\cong G/K$ be a symmetric space of non-compact type,~and $\g{g}=\g{k}\oplus\g{p}$ the Cartan decomposition determined by the point $o\in M$ whose isotropy subgroup of $G$ is precisely $K$. Given a maximal abelian subspace $\g{a}$ of $\g{p}$, and a covector $\lambda\in \g{a}^*$, we define $\g{g}_\lambda=\{X\in\g{g}:[H,X]=\lambda(H)X \text{ for all }H\in\g{a}\}$. The choice of $\g{a}$~determines the set of restricted roots, $\Delta=\{\lambda\in\g{a}^*:\g{g}_\lambda\neq 0\}\setminus\{0\}$, and the restricted root space~decomposition of $\g{g}$, namely $\g{g}=\g{g}_0\oplus(\bigoplus_{\lambda\in\Delta}\g{g}_\lambda)$. The relation $[\g{g}_\lambda,\g{g}_\mu]\subset\g{g}_{\lambda+\mu}$ holds for any $\lambda$, $\mu\in\g{a}^*$. 

When $\g{a}^*$ is endowed with the inner product induced by the Killing form of $\g{g}$, the subset $\Delta$ of $\g{a}^*$ turns out to be a (possibly non-reduced) root system. Let $\Delta^+$ be a set of positive roots for  $\Delta$, defined as the subset of roots lying at one side of a linear hyperplane in $\g{a}^*$ not intersecting $\Delta$. Thus, $\Delta$ is the disjoint union of $\Delta^+$ and $-\Delta^+$. We denote by $\Pi\subset\Delta^+$ the corresponding set of simple roots, which is made of all positive roots that are not sums of two positive roots. In particular, $\Pi$ is a basis of $\g{a}^*$, and every $\lambda\in\Delta$ is a linear combination of elements in $\Pi$ whose coefficients are all either non-negative or non-positive integers.

We consider the nilpotent Lie algebra $\g{n}=\bigoplus_{\lambda\in\Delta^+}\g{g}_\lambda$, and the connected Lie subgroups $A$, $N$ and $AN$ of $G$ with Lie algebras $\g{a}$, $\g{n}$ and $\g{a}\oplus\g{n}$, respectively. The Iwasawa decomposition theorem for $G$ implies that $M\cong G/K$ is diffeomorphic to the solvable Lie group $AN$. We denote by $\langle\cdot, \cdot \rangle$ both the Riemannian metric on $M$ and the pullback metric on $AN$ that makes $M$ and $AN$ (as well as $T_oM$ and $\g{a}\oplus\g{n}$) isometric. It turns out that such metric on $AN$ is left-invariant. 
Thus, by Koszul's formula, the Levi-Civita connection of $AN$ is determined by the expression $2\langle \nabla_X Y,Z\rangle =\langle [X,Y],Z\rangle + \langle [Z,X],Y\rangle +\langle X,[Z,Y]\rangle$, for any $X$, $Y$, $Z\in\g{a}\oplus\g{n}$. In particular, if $H\in\g{a}$ and $X\in\g{a}\oplus\g{n}$, we have 
\begin{equation}\label{equation:levi:civita}
\nabla_X H=-[H,X].
\end{equation}

For each $\lambda\in\g{a}^*$, we define the vector $H_\lambda\in\g{a}$ by the relation $\langle H_\lambda,H\rangle=\lambda(H)$ for all $H\in\g{a}$. Particular instances of this definition are the vectors $H_\lambda\in\g{a}$ where $\lambda\in\Delta\subset\g{a}^*$ is a root, and also the vector $H_\delta\in\g{a}$, where
$
\delta=\frac{1}{2}\sum_{\lambda\in\Delta^+}(\dim \g{g}_\lambda) \lambda\in\g{a}^*$. We will denote also by $\langle \cdot,\cdot\rangle$ the inner product on $\g{a}^*$ given by $\langle \lambda,\mu\rangle=\langle H_\lambda, H_\mu\rangle$, for any $\lambda,\mu\in\g{a}^*$. With respect to this inner product of $\g{a}^*$, $\Delta$ is also a (possibly non-reduced) root system on $\g{a}^*$.

If $k$ is an isometry of $M$ fixing $o$ and such that $\Ad(k)\g{a}=\g{a}$, 
then it is a standard fact that $\Ad(k)$ permutes the root spaces $\g{g}_\lambda$ with $\lambda\neq 0$, and given $\lambda$, $\mu\in \Delta$, the relations
\begin{equation}\label{eq:auto}
\Ad(k)H_\lambda=H_\mu, \qquad \mu=\lambda\circ\Ad(k^{-1})\vert_\g{a}, \qquad \Ad(k)\g{g}_\lambda=\g{g}_\mu
\end{equation}
are equivalent to each other. This follows from the definitions and the fact that $\Ad(k)$ is an automorphism of $\g{g}$ and a linear isometry when restricted to $\g{a}$.

\medskip
In the following subsections we will prove Theorems~\ref{th:explicit} and~\ref{th:extension}. As already mentioned, Theorem~\ref{th:first} follows directly from Theorem~\ref{th:explicit}. 

Let $S$ be the connected Lie subgroup of $AN$ with Lie algebra $\g{s}=\g{b}\oplus\g{n}$, where $\g{b}$ is a  proper vector subspace of $\g{a}$. We put $\g{b}^\perp:=\g{a}\ominus\g{b}$. Hereafter, $\ominus$ denotes orthogonal complement.

\subsection{Extrinsic geometry and non-austerity of the focal set}\label{subsec:extrinsic}
In this subsection we will prove Theorem~\ref{th:explicit}~(i).
We first calculate the shape operator $\cal{S}$ of $S\cdot o$ as a submanifold of $M$. By homogeneity, it will be enough to calculate $\cal{S}$ at $o$. Thanks to the isometry $M\cong AN$ and using~\eqref{equation:levi:civita}, $\cal{S}$ can be computed as $\cal{S}_\xi X=-(\nabla_X \xi)^\top = -[\xi, X]^\top$, where $\xi\in \g{b}^\perp$ and $X\in\g{s}$ are left-invariant vector fields on $AN$, and $^\top$ denotes orthogonal projection onto the tangent space $\g{s}$. Using~\eqref{equation:levi:civita}, we have
\begin{equation}\label{eq:shape}
	\cal{S}_\xi H=0 \quad \text{and}\quad \cal{S}_\xi X_\lambda=\lambda(\xi)X_\lambda, \qquad \text{for any } \xi\in\g{b}^\perp,\,H\in\g{b} \text{ and } X_\lambda\in\g{g}_\lambda\subset \g{n}.
\end{equation}
Hence, the mean curvature vector $\cal{H}$ of $S\cdot o$ is determined by the relation $\langle \cal{H},\xi\rangle=\tr \cal{S}_\xi=\sum_{\lambda\in \Delta^+} \dim\g{g}_\lambda\lambda(\xi)=2\langle H_\delta,\xi\rangle$, for each $\xi\in \g{b}^\perp$. Thus, $S\cdot o$ is minimal if and only if $H_\delta\in\g{s}$, or equivalently $H_\delta\in\g{b}$.  This proves the first claim in Theorem~\ref{th:explicit}~(i).

\begin{remark}\label{rem:Einstein}
	The condition $H_\delta\in\g{b}\subset\g{a}$ also arises in the study of Einstein solvmanifolds. Indeed, this assumption means that $\g{s}$ contains the mean curvature vector $\cal{H}=2H_\delta$ of the solvmanifold $AN$, and thus, if the symmetric metric on $M$ is Einstein (for example, if $M$ is irreducible), then the orbit $S\cdot o$ is known to be Einstein as well, see~\cite[Theorem~4.18]{Heber}. Moreover, in this situation, the $S$-action on $M$ is free and polar with mutually isometric minimal orbits (see~\S\ref{subsec:extension}). The proof of the  existence of such an action on any Einstein manifold of negative scalar curvature in the presence of symmetry is one of the key steps in the recent outstanding proof of the Alekseevskii conjecture by B\"ohm and Lafuente~\cite{BL:arxiv}.
\end{remark}

From now on in this subsection and the next one (except in Remark~\ref{rem:rank2} below), we will assume that $H_\delta\in\g{b}$ and $\dim \g{b}^\perp\geq 2$, and hence $\rank M\geq 3$. 

\medskip
We will now prove the second claim in Theorem~\ref{th:explicit}~(i) by showing that $S\cdot o$ is not austere if $\g{b}=\R H_\delta$, and that $S\cdot o$ is not austere for a generic choice of a subspace $\g{b}$ of $\g{a}$ containing $H_\delta$. From~\eqref{eq:shape} we know that the principal curvatures of $S\cdot o$, regarded as linear functionals on the normal space $\g{b}^\perp$, are $0$ and the restrictions $\lambda\vert_{\g{b}^\perp}$ of the positive roots $\lambda\in \Delta^+$. Let $\lambda_1,\dots, \lambda_\ell$ be an enumeration of the positive roots in $\Delta^+$, where each root $\lambda_i$ is repeated as many times as the dimension $\dim\g{g}_{\lambda_i}$ of its associated root space; in particular, $\ell=\dim\g{n}$. 
Thus, $S\cdot o$ is austere if and only if
\begin{equation}\label{eq:austere_characterization}
\text{there is a permutation } \sigma \text{ of } \lambda_1,\dots, \lambda_\ell \text{ with } \left(\lambda_i+\sigma(\lambda_i)\right) \vert_{\g{b}^\perp}=0 \text{ for all } i\in\{1,\dots,\ell\}. 
\end{equation}
	
Note that if $\lambda$, $\mu\in\Delta^+$, then $\lambda+\mu$ vanishes on $\g{b}^\perp$ if and only if $H_{\lambda+\mu}\in\g{b}$. Therefore, if for some $\lambda\in\Delta^+$ we have $H_{\lambda+\mu}\notin\g{b}$ for all $\mu\in\Delta^+$, then $S\cdot o$ cannot be austere. Assume for a moment that 
\begin{equation}\label{eq:collinear}
\text{there is }\lambda\in\Delta^+\text{ such that }\lambda+\mu \text{ is not collinear to } \delta \text{ for any } \mu\in\Delta^+.
\end{equation}
Then it follows that $S\cdot o$ is not austere for the choice $\g{b}=\R H_\delta$. More generally, any subspace $\g{b}$ of $\g{a}$ such that
\begin{equation}\label{eq:sufficient_non_austere}
	H_\delta\in\g{b}, \qquad \text{and}\qquad H_{\lambda+\mu}\notin\g{b} \text{ for any }\mu\in\Delta^+
\end{equation}
will produce a non-austere minimal $S\cdot o$, where $\lambda$ is a fixed root satisfying~\eqref{eq:collinear}. 
Observe that the $k$-dimensional subspaces $\g{b}$ of $\g{a}$ satisfying the conditions in~\eqref{eq:sufficient_non_austere} constitute an open and dense subset of the space of $k$-dimensional subspaces of $\g{a}$ containing $H_\delta$, for any $k\in\{1,\dots, \rank M-1\}$. This is because, for any given $\lambda\in\Delta^+$ as in~\eqref{eq:collinear}, $\{H_{\lambda+\mu}: \mu\in\Delta^+\}$ is a finite subset of $\g{a}$ not intersecting $\R H_\delta$.

Thus, in order to conclude the proof of Theorem~\ref{th:explicit}~(i), we just have to show that \eqref{eq:collinear} is true for any (possibly non-reduced) root system of rank at least $3$. This is probably well known, but since we did not find an appropriate reference in the literature, we shall give a proof here. We will make frequent use of known properties of root systems, particularly those in~\cite[Proposition~2.48]{Knapp}. We will also use the fact that $\langle \delta, \nu\rangle >0$ for all $\nu\in\Delta^+$. This can be easily shown by adapting the proof of~\cite[Proposition~2.69]{Knapp}. Thus, \eqref{eq:collinear} will be proved if we show that
\begin{equation}\label{eq:leq}
	\text{there is } \lambda\in\Delta^+ \text{ such that for each } \mu\in\Delta^+ \text{ we have } \langle \lambda+\mu,\nu\rangle\leq 0 \text{ for some }\nu\in\Delta^+.
\end{equation}

We first consider the case that $\Delta=\Delta_1\oplus\Delta_2$ is a reducible root system. Since $\Delta$ has rank at least $3$, we can assume that the (possibly reducible) factor $\Delta_1$ has rank $2$ or higher. Fix $\lambda\in\Pi\cap \Delta_1$ a simple root in the factor $\Delta_1$, and let $\mu\in\Delta^+$ be arbitrary. If $\mu\in\Delta_1$, then any $\nu\in\Delta_2^+$ satisfies $\langle \lambda+\mu,\nu\rangle=0$. Otherwise, if $\mu\in\Delta_2$, since $\rank \Delta_1\geq 2$, there exists $\nu\in \Pi\cap\Delta_1$, $\nu\neq \lambda$, and hence $\langle \lambda +\mu,\nu\rangle=\langle\lambda,\nu\rangle\leq 0$, since $\lambda$ and $\nu$ are simple. This shows~\eqref{eq:leq}, and hence~\eqref{eq:collinear}, in the reducible case.

Now assume that $\Delta$ is irreducible. 
Let $\lambda\in\Delta^+$ be such that $|\lambda|\geq |\mu|$ for every $\mu\in \Delta^+$, and $\lambda\in \Pi$ or $\lambda/2\in \Pi$ (if $\Delta$ is non-reduced). Let $\mu\in\Delta^+$ be arbitrary. 
If $\mu\in \Pi$ or if $\mu/2\in \Pi$ (in case $\mu$ is non-reduced), then any simple root $\nu\in \Pi\setminus\spann\{\lambda,\mu\}$ (which exists since $\rank\Delta\geq 3$) satisfies $\langle \lambda,\nu\rangle\leq 0$ and $\langle \mu,\nu\rangle \leq 0$, from where we get~\eqref{eq:leq}. Otherwise, assume that $\R\mu\cap\Pi=\emptyset$. Take $\nu\in \Pi$ joined to $\lambda$ (or to $\lambda/2$ if $\lambda$ is non-reduced) in the Dynkin diagram associated with $\Pi$. Then $\langle \lambda,\nu\rangle=-|\lambda|^2/2$, and $\langle \mu,\nu\rangle\leq \max\{|\mu|^2,|\nu|^2\}/2\leq |\lambda|^2/2$, where we have used that neither $\lambda$ nor $\mu$ are proportional to $\nu$.  This implies~\eqref{eq:leq}, and thus we have concluded the proof of~\eqref{eq:collinear}, and hence of Theorem~\ref{th:explicit}~(i).

\begin{remark}\label{rem:rank2}
	The assumption that $\rank M\geq 3$ is essential for the proof of the non-austerity claim in Theorem~\ref{th:explicit}~(i). Specifically, it was crucial in the proof of~\eqref{eq:leq}. We did not consider the case $\rank M = 2$ (with $\dim \g{b}=1$) in Theorem~B, as it would lead to a minimal leaf $S \cdot o$ of a homogeneous (hence isoparametric) codimension one regular foliation, and these were classified in~\cite{BDT:jdg}. However, in this case, the austerity of the hypersurface $S\cdot o$ depends on $M$. Indeed, using the characterization in \eqref{eq:austere_characterization}, one can check for instance that if $M$ is of type $A_2$, then $S\cdot o$ is austere, whereas if $M$ is of type $B_2$, then it is not.
\end{remark}

\begin{remark}\label{rem:austere}
	The problem of completely determining the subspaces $\g{b}$ of $\g{a}$ that produce austere orbits $S\cdot o$ seems to be a non-straightforward combinatorial problem. In this article we content ourselves with proving that non-austere examples are generic, and with giving some examples (necessarily with $\R H_\delta\subsetneq\g{b}$) which show that austere orbits $S\cdot o$ indeed~exist. 
	
	Consider $M=SL_5/SO_5$, let $\Delta^+$ be a set of positive restricted roots (which has $10$ elements), and $\Pi=\{\alpha_1,\alpha_2,\alpha_3,\alpha_4\}\subset\Delta^+$ the corresponding set of simple roots (ordered in the standard way, so that $\alpha_1$ and $\alpha_4$ correspond to the extremal nodes of the $A_4$-Dynkin diagram). Then $\g{b}=\spann\{H_\lambda, H_{\alpha_2+\alpha_3}\}$ yields an austere orbit $S\cdot o$, where $\lambda=\sum_{i=1}^4\alpha_i$ is the highest root. Indeed the permutation $\sigma$ of $\Delta^+$ given in cycle notation by
	$(\alpha_1\;\lambda-\alpha_1)(\alpha_2\;\alpha_3)(\alpha_4\;\lambda-\alpha_4)(\alpha_1+\alpha_2\;\alpha_3+\alpha_4)$ satisfies the characterization in~\eqref{eq:austere_characterization}.
	For an example in the reducible setting, let $M=(\R H^2)^4$ and denote by $\alpha_1,\alpha_2,\alpha_3,\alpha_4$ the roots in $\Delta^+=\Pi$. Then, the subspace $\g{b}=\spann\{H_{\alpha_1+\alpha_2},H_{\alpha_3+\alpha_4}\}$ of $\g{a}$ produces an austere orbit $S\cdot o$, as the permutation $\sigma$ of $\Delta^+$ given by $(\alpha_1\; \alpha_2)(\alpha_3\;\alpha_4)$ satisfies~\eqref{eq:austere_characterization}.
	
\end{remark}

\begin{remark}\label{rem:CPC}
	An important subclass of austere submanifolds is that determined by the property that all shape operators for unit normal vectors are isospectral. These were called CPC submanifolds in~\cite{BS:cag@gmail.com}. Their importance stems from the fact that focal sets of isoparametric families of hypersurfaces with constant principal curvatures have this property~\cite{GT:asian}, and indeed, in the context of spaces of constant curvature,  the classification of CPC submanifolds is equivalent to that of isoparametric hypersurfaces. One may ask if, for some choices $\g{b}\subset \g{a}$, our construction method produces CPC orbits $S\cdot o$. The argument in the last paragraph of \S\ref{subsec:isoparametric} will show that this is never the case.
\end{remark}

\subsection{Tubes around $S\cdot o$ are isoparametric with non-constant principal curvatures}\label{subsec:isoparametric}
In this subsection we will prove Theorem~\ref{th:explicit}~(ii). To this purpose, we will calculate the extrinsic geometry of the distance tubes around $S\cdot o$ using Jacobi field theory; we refer the reader to~\cite[\S10.2]{BCO} for details on this method. Let $\gamma$ be a unit speed geodesic in $M$ with $\gamma(0)=o$ and $\dot{\gamma}(0)=\xi\in \nu_o (S\cdot o)$, where $\nu (S\cdot o)$ is the normal bundle of $S\cdot o$. Given $X\in T_o M$, we will denote by $\cal{P}_X$ the parallel vector field along $\gamma$ with $\cal{P}_X(0)=X$. 
For each $X\in T_oM\cong \g{a}\oplus\g{n}$ orthogonal to $\xi\in \nu_o(S\cdot o)\cong\g{b}^\perp\subset\g{a}$, let $J_X$ be the adapted Jacobi vector field along $\gamma$ that is solution to the initial value problem
\[
J_X''+R(J_X,\dot{\gamma})\dot{\gamma}=0, \quad J_X(0)=X^\top, \quad J_X'(0)=X^\perp-\cal{S}_{\xi} X^\top,
\]
where $R$ is the curvature tensor of $M$, and $^\top$, $^\perp$ denote the projections onto the tangent and normal spaces of $S\cdot o$, respectively.

Let $H\in\g{a}$ and $X_\lambda\in\g{g}_\lambda$, $\lambda\in \Delta^+$, be arbitrary vectors, which we regard as tangent vectors to $M$ at $o$. It is well known (see~\cite[\S2.15]{eberlein}, and note that $\xi\in\g{a}$) that the parallel vector fields $\cal{P}_H$ and $\cal{P}_{X_\lambda}$ along the homogeneous geodesic $\gamma$ are eigenvectors of the Jacobi operator $R_{\dot{\gamma}}=R(\cdot, \dot{\gamma})\dot{\gamma}$ with respective eigenvalues $0$ and $-\lambda(\xi)^2$. Thus, given arbitrary vectors $\eta\in \g{b}^\perp\ominus\R\xi$, $H\in\g{b}\subset\g{s}$ and $X_\lambda\in\g{g}_\lambda\subset \g{s}$, one can easily verify, taking \eqref{eq:shape} into account, that the adapted Jacobi fields $J_\eta$, $J_H$ and $J_{X_\lambda}$ are given~by 
\begin{equation}\label{eq:jacobi_fields}
J_\eta(t)=t \,\cal{P}_\eta(t),\quad J_H(t)=\cal{P}_H(t), \quad J_{X_\lambda}(t)=
e^{-t\lambda(\xi)}\cal{P}_{X_\lambda}(t).
\end{equation}
Since for any $t>0$ and unit normal vector $\xi\in\nu_o(S\cdot o)$, the set $\{J_X(t):X\in {T_oM\ominus\R\xi}\}$ spans a $(\dim M-1)$-dimensional subspace of $T_{\gamma(t)}M$, and taking into account the homogeneity of $S\cdot o$, it follows that the tube $W^t:=\{\exp (t\xi):\xi\in\nu(S\cdot o), |\xi|=1\}$ of radius $t$ around $S\cdot o$ is indeed a hypersurface of $M$, for any $t>0$. Moreover, the shape operator $\cal{S}^t$ of $W^t$ with respect to the normal vector $\dot{\gamma}(t)$ can be calculated by $\cal{S}^tJ_X(t)=-J_X'(t)$, for each $X\in T_oM\ominus\R\xi$. Therefore, inserting~\eqref{eq:jacobi_fields} into the previous relation, we obtain
\begin{equation}\label{eq:shape^t}
\cal{S}^t\cal{P}_\eta(t)=-\frac{1}{t}\cal{P}_\eta(t), \qquad 
\cal{S}^t\cal{P}_H(t)=0, \qquad 
\cal{S}^t\cal{P}_{X_\lambda}(t)=\lambda(\xi)\cal{P}_{X_\lambda}(t),
\end{equation}
for any $\eta\in \g{b}^\perp\ominus\R \xi$, $H\in\g{b}\subset\g{s}$ and $X_\lambda\in\g{g}_\lambda\subset \g{s}$. Thus, the mean curvature of $W^t$ is 
\[
\cal{H}^t=\tr \cal{S}^t=-\frac{1}{t}(\dim \g{b}^\perp-1)+\sum_{\lambda\in\Delta^+}\dim\g{g}_\lambda\lambda(\xi)=-\frac{1}{t}(\dim \g{b}^\perp-1),
\]
where in the last equality we have used the assumption that $H_\delta$ is orthogonal to any normal vector $\xi$. This shows that each tube $W^t$ has constant mean curvature. As discussed in the introduction, each $W^t$ is the union of the $S$-orbits passing through the points of a sphere of radius $t$  centered at $o$ in the Euclidean space $\Sigma:=\Exp(\g{b}^\perp)\cdot o$. Since the action of $S\subset AN$ on $M$ is free, and by construction the tubes $W^t$ are equidistant to each other and to $S\cdot o$, we conclude that $\{W^t:t>0\}\cup\{S\cdot o\}$ is a isoparametric family of hypersurfaces with focal set $S\cdot o$.

Finally, the non-constancy of the principal curvatures of $W^t$, and hence the inhomogeneity of $W^t$, follows directly from~\cite[Theorem~1.2]{GT:asian} when $S \cdot o$ is not austere. However, we can provide an easy argument, based on \eqref{eq:shape^t}, which works in general, independently of whether $S \cdot o$ is austere or not. Indeed, if the principal curvatures of some tube $W^t$ around $S \cdot o$ were constant, the set $\{\lambda(\xi):\lambda\in\Delta^+,\xi\in\g{b}^\perp, |\xi|=1\}$ would be finite. Note that such set is also constant if $S\cdot o$ is assumed to be CPC (cf.~Remark~\ref{rem:CPC}). Since $\dim\g{b}^\perp\geq 2$, this would necessarily imply $\g{b}^\perp\subset \ker \lambda$ for all $\lambda\in\Delta^+$, which is impossible, as $\Delta^+$ generates $\g{a}^*$. This concludes the proof of item (ii) of Theorem~\ref{th:explicit} and of Remark~\ref{rem:CPC}.

\subsection{Congruence of the examples}\label{subsec:congruence}
This subsection contains the proof of Theorem~\ref{th:explicit}~(iii). More specifically, let $\g{b}_1$, $\g{b}_2$ be two vector subspaces of $\g{a}$, and let $S_i$ be the connected Lie subgroup of $AN$ with Lie algebra $\g{s}_i=\g{b}_i\oplus\g{n}$, $i\in\{1,2\}$. We will prove that if the orbits $S_1\cdot o$ and $S_2\cdot o$ are congruent, then $\g{b}_2=\varphi(\g{b}_1)$ for some linear isometry $\varphi$ of $\g{a}$ of the form $\varphi=\Ad(k)\vert_\g{a}$, $k\in N_{\widehat{K}}(\g{a})$. Here $N_{\widehat{K}}(\g{a})$ denotes the normalizer of $\g{a}$ in the full isotropy group $\widehat{K}=\{k\in \mathrm{Isom}(M):k(o)=o\}$. Note that the group $\{\Ad(k)\vert_\g{a}:k\in N_{\widehat{K}}(\g{a})\}$ leaves the finite set $\{H_\lambda\in\g{a}:\lambda\in \Delta\}$ invariant, see~\eqref{eq:auto}. It is therefore a finite group, as the $H_\lambda$, $\lambda\in \Delta$, generate $\g{a}$.

We will not assume any restriction on the rank of $M$ or on the subspaces $\g{b}_1$ and $\g{b}_2$. Observe that Theorem~\ref{th:explicit}~(iii) will then follow directly, since two choices of subspaces $\g{b}_1$, $\g{b}_2$ of $\g{a}$ in the conditions of Theorem~\ref{th:explicit} produce congruent isoparametric families if and only if the orbits $S_1\cdot o$ and $S_2\cdot o$ are congruent. 

Assume that $S_1\cdot o$ and $S_2\cdot o$ are congruent submanifolds of $M$. By composing with an element of $S_2$ if necessary, we can assume that there is an isometry $k\in\widehat{K}$ of $M$ fixing $o$ and such that $k(S_1\cdot o)=S_2\cdot o$. Then $k_*(\nu_o(S_1\cdot o))=\nu_o(S_2\cdot o)$, which translates into $\Ad(k)\g{b}_1^\perp=\g{b}_2^\perp$, as $\g{b}_i^\perp=\g{a}\ominus\g{b}_i$ is a subspace of $\g{p}$, for each $i\in\{1,2\}$. Now, Theorem~\ref{th:explicit}~(iii) follows directly from the following

\begin{lemma}\label{lemma:congruence}
Let $\g{b}_1$, $\g{b}_2$ be subspaces of $\g{a}$ such that $\Ad (k) \g{b}_1^\perp=\g{b}_2^\perp$ for some $k\in\widehat{K}$. Then, there exists $\tilde{k} \in N_{\widehat{K}}(\g{a})$ such that $\Ad(\tilde{k})\g{b}_1=\g{b}_2$.
\end{lemma}

\begin{proof}
For each $i\in\{1,2\}$, define $\Delta_i^+:=
\{\lambda\in \Delta^+:\lambda({\g{b}_i^\perp})=0\}$ and $\g{a}_i:=\bigcap_{\lambda\in \Delta_i^+}\ker\lambda$. By definition, $\g{b}_i^\perp\subset\g{a}_i$. It turns out that the centralizer of $\g{b}_i^\perp$ in $\g{p}$ is a Lie triple system of the form $Z_\g{p}(\g{b}_i^\perp)=\g{a}_i\oplus\g{c}_i$, where $\g{c}_i=(\g{a}\ominus\g{a}_i)\oplus\bigl(\bigoplus_{\lambda\in\Delta_i^+}\g{p}_\lambda\bigr)$, $\g{p}_\lambda=(1-\theta)\g{g}_\lambda\subset\g{p}$, and $\theta$ is the Cartan involution associated with the decomposition $\g{g}=\g{k}\oplus\g{p}$; see~\cite[\S2.20]{eberlein} for details. Both $\g{a}_i$ and $\g{c}_i$ are also (mutually orthogonal) Lie triple systems, corresponding to a  totally geodesic flat submanifold $\Exp(\g{a}_i)\cdot o\cong\R^{r-r_i}$ and a totally geodesic semisimple symmetric space of non-compact type $C_i:=\Exp(\g{c}_i)\cdot o$ of rank $r_i$, respectively, where $r:=\rank M=|\Pi|$ and $r_i=\dim(\g{a}\ominus\g{a}_i)$ is the maximum number of linearly independent roots in $\Delta_i^+$. Since $\Ad(k)\g{b}_1^\perp=\g{b}_2^\perp$, we also have $\Ad(k)Z_\g{p}(\g{b}_1^\perp)=Z_\g{p}(\g{b}_2^\perp)$, and hence $\Ad(k)(\g{a}_1\oplus\g{c}_1)=\g{a}_2\oplus\g{c}_2$. Thus, $k$ maps the totally geodesic submanifold $\R^{r-r_1}\times C_1$ to $\R^{r-r_2}\times C_2$, and then $k$ necessarily preserves the Euclidean and semisimple factors, from where $\Ad(k)\g{a}_1=\g{a}_2$ and $\Ad(k)\g{c}_1=\g{c}_2$.

Let $K_i$ be the connected component of the identity of the isotropy group of $C_i$ at $o$. Its Lie algebra is $\g{k}_i=[\g{c}_i,\g{c}_i]$. Note that $\g{a}\ominus\g{a}_2$ and $\Ad(k)(\g{a}\ominus\g{a}_1)\subset\Ad(k)\g{c}_1=\g{c}_2$ are maximal abelian subspaces of $\g{c}_2$. Since any two maximal abelian subspaces of the tangent space of a symmetric space are conjugate, there exists $k'\in K_2$ such that the element $\tilde{k}:=k'k\in \widehat{K}$ satisfies $\Ad(\tilde{k})(\g{a}\ominus\g{a}_1)=\g{a}\ominus\g{a}_2$. As $\g{k}_2=[\g{c}_2,\g{c}_2]$ centralizes $\g{a}_2$, then $\Ad(k')$ leaves $\g{a}_2$ (and hence $\g{b}_2^\perp$) pointwise invariant. This, together with $\Ad(k)\g{b}_1^\perp=\g{b}_2^\perp$, $\Ad(k)\g{a}_1=\g{a}_2$ and $\Ad(\tilde{k})(\g{a}\ominus\g{a}_1)=\g{a}\ominus\g{a}_2$, implies $\Ad(\tilde{k})\g{a}=\g{a}$ and $\Ad(\tilde{k})\g{b}_1=\g{b}_2$, as we wanted to show.
\end{proof}

\subsection{Extension method and proof of Theorem~\ref{th:extension}}\label{subsec:extension}
The Lie group $AN$ acts simply transitively on $M$. Hence, the connected Lie subgroup $S$ of $AN$ with Lie algebra $\g{s}=\g{b}\oplus\g{n}$ acts freely on $M$. Moreover, the $S$-action on $M$ is hyperpolar with section $\Sigma:={\Exp(\g{b}^\perp)\cdot o}\cong \R^d$, where $d=\dim\g{b}^\perp=\rank M-\dim\g{b}$. Indeed, such an $S$-action is one of the examples of hyperpolar actions with no singular orbits classified in~\cite{BDT:jdg}. Furthermore, $\Sigma$ intersects each $S$-orbit exactly once, as follows from the fact that $AN= \Exp(\g{b}^\perp) S=S\Exp(\g{b}^\perp)$ acts simply transitively on $M$. Since $H_\delta\in \g{b}$, the orbit $S\cdot o$ is a minimal submanifold of $M$, as shown in~\S\ref{subsec:extrinsic}. By polarity, any $S$-orbit is of the form $S\cdot a(o)=a(a^{-1}Sa\cdot o)=a(S\cdot o)$, for some $a\in\Exp(\g{b}^\perp)$, where we have used that the group $A$ normalizes $S$. Thus, all $S$-orbits are mutually congruent, and hence minimal (cf.~\cite{BDT:jdg}). 
Therefore, the action of $S$ on $M$ satisfies the hypotheses of \cite[Theorem~2.1]{DV:imrn}, from where the first claim in item (i) as well as item (ii) in Theorem~\ref{th:extension} follow immediately (for item (ii), note that polar singular Riemannian foliations on Euclidean spaces are hyperpolar and isoparametric).

As discussed in~\cite{DV:imrn}, the shape operator of an extended submanifold $S\cdot P$ of $M$, where $P$ is a submanifold of $\Sigma$, is block diagonal, with one block corresponding to the shape operator of $P$, and the other block to the shape operator of the orbit $S\cdot o$. Since the latter never vanishes identically, as follows from~\eqref{eq:shape} and the argument in the last paragraph of~\S\ref{subsec:isoparametric}, we have that $S\cdot P$ is never totally geodesic, which completes the proof of Theorem~\ref{th:extension}~(i).

We will finally prove item (iii) of Theorem~\ref{th:extension}. Let $P$ be a submanifold of $\Sigma$, and $S\cdot P$ the extended submanifold of $M$. We will show that $S\cdot P$ is a homogeneous submanifold of $M$ if and only if $P$ is an affine subspace of $\Sigma\cong\R^d$; the analogous claim for foliations stated in item (iii) follows immediately. Note that $P$ is an affine subspace of $\Sigma=\Exp(\g{b}^\perp)\cdot o$ if and only if $P$ is an orbit $L\cdot p$, $p\in \Sigma$, of the connected Lie subgroup $L$ of $A$ whose Lie algebra is a vector subspace $\g{l}$ of $\g{b}^\perp$. In this case, $S\cdot P$ is the orbit of the group $LS$ through $p$, and thus it is homogeneous. 

Conversely, assume that $S\cdot P$ is homogeneous. By applying an isometry in $A$ if necessary, we can assume that $o\in P$. Thus, in particular, for each $p\in P=(S\cdot P)\cap \Sigma$ there exists $g_p\in \mathrm{Isom}(M)$ such that $g_p(o)=p$ and $g_p(S\cdot P)=S\cdot P$. Moreover, for each such $p$ there is exactly one $a_p\in \Exp(\g{b}^\perp)\subset A$ such that $a_p(p)=o$. Recall that $\widehat{K}=\{k\in \mathrm{Isom}(M):k(o)=o\}$. Define $k:=a_pg_p\in\widehat{K}$, which satisfies $k(S\cdot P)=a_p(S\cdot P)= S\cdot a_p(P)$, since $a_p\in A$ normalizes $S$. Consider the normal spaces at $o$ to the submanifolds $S\cdot P$ and $S\cdot a_p(P)$. Since such normal spaces are tangent to $\Sigma$, we can identify $\nu_o(S\cdot P)\cong \g{b}_1^\perp$ and $\nu_o(S\cdot a_p(P))\cong \g{b}_2^\perp$, for subspaces $\g{b}_1^\perp$ and $\g{b}_2^\perp$ of $\g{b}^\perp\subset\g{a}$. Since $k_*(\nu_o(S\cdot P))=\nu_o(S\cdot a_p(P))$, we have $\Ad(k)\g{b}_1^\perp=\g{b}_2^\perp$. Thus, from Lemma~\ref{lemma:congruence} we deduce the existence of an element $\tilde{k}\in N_{\widehat{K}}(\g{a})$ such that  $\Ad(\tilde{k})\g{b}_1=\g{b}_2$, and hence $\Ad(\tilde{k})\g{b}_1^\perp=\g{b}_2^\perp$. Since $\{\Ad(\hat{k})\vert_\g{a}:\hat{k}\in N_{\widehat{K}}(\g{a})\}$ is a finite group of linear isometries of $\g{a}$, then there is only a finite number of normal spaces $\nu_o(S\cdot a_p(P))$ for $p\in P$. But by the connectedness of $P$ and the continuous dependence of $a_p$ on $p$, we must have $\nu_o(S\cdot a_p(P))=\nu_o(S\cdot P)$, for all $p\in P$. Equivalently, regarding $P$ and $a_p(P)$ as submanifolds of $\Sigma\cong\R^d$, we have $\nu_o(a_p(P))=\nu_o(P)$, for all $p\in P$. But the only submanifolds of a Euclidean space $\R^d$ with this property (where the isometries $a_p$ are translations of $\R^d$) are open parts of affine subspaces. The homogeneity of $S\cdot P$ implies that $P$ must indeed be a complete affine subspace. This concludes the proof of Theorem~\ref{th:extension}.

\enlargethispage{2\baselineskip}

\end{document}